\documentclass[12pt]{article}
\usepackage{graphicx}
\usepackage{amsfonts}
\usepackage{latexsym}
\usepackage{amsmath}
\usepackage{fancyhdr}

\makeatletter
\@addtoreset{figure}{section}
\def\thefigure{\thesection.\@arabic\c@figure}
\@addtoreset{table}{section}
\def\thetable{\thesection.\@arabic\c@table}

\def\@sect#1#2#3#4#5#6[#7]#8{\ifnum #2>\c@secnumdepth
     \def\@svsec{}\else
     \refstepcounter{#1}\edef\@svsec{\csname the#1\endcsname.\hskip .75em
}\fi
     \@tempskipa #5\relax
      \ifdim \@tempskipa>\z@
        \begingroup #6\relax
          \@hangfrom{\hskip #3\relax\@svsec}{\interlinepenalty \@M #8\par}%
        \endgroup
       \csname #1mark\endcsname{#7}\addcontentsline
         {toc}{#1}{\ifnum #2>\c@secnumdepth \else
                      \protect\numberline{\csname the#1\endcsname}\fi
                    #7}\else
        \def\@svsechd{#6\hskip #3\@svsec #8\csname #1mark\endcsname
                      {#7}\addcontentsline
                           {toc}{#1}{\ifnum #2>\c@secnumdepth \else
                             \protect\numberline{\csname the#1\endcsname}\fi
                       #7}}\fi
     \@xsect{#5}}
\def\@begintheorem#1#2{\it \trivlist \item[\hskip \labelsep{\bf #1\ #2.}]}
\def\section{\@startsection {section}{1}{\z@}{-3.5ex plus -1ex minus
 -.2ex}{2.3ex plus .2ex}{\normalsize\bf}}

\begin{document}

\title{The Structure of The Group of Polynomial Matrices\\
Unitary In The Indefinite Metric of Index 1}
%uncomment to remove date
\date{} 
\maketitle

\begin{center}
%\small
\author{Boris D. Lubachevsky\\
${ }$\\
{\em bdl@bell-labs.com}\\
${ }$\\
Bell Laboratories \\
600 Mountain Avenue \\
Murray Hill, New Jersey}
\end{center}

\setlength{\baselineskip}{0.995\baselineskip}
\normalsize
\vspace{0.5\baselineskip}
\vspace{1.5\baselineskip}
%\end{center}

% ----------------------------------------------------------------
\begin{abstract}
We consider the group $\cal M$ of all polynomial $\nu \times \nu$ matrices 
$U(\omega) = \sum_{\iota=0}^{\kappa} U_{\iota}\, \omega^{\iota}$, 
$\kappa=0, 1, 2,...$,
that satisfy equation 
$U(\omega)\cdot D\cdot U(\omega)^*=D$ with diagonal
$\nu \times \nu$
matrix $D = \mbox{diag}\{-1,1,1,...1\}$. 
Here $\nu \ge 2$,
$U(\omega)^*
 = \sum_{\iota=0}^{\kappa} U_\iota^* \omega^\iota$ and $U_\iota^*$ 
is the Hermitean conjugate of $U_\iota$.
We show that the subgroup 
${\cal M}_0$ of those $U(\omega) \in \cal M$,
that are normalized by the condition
$U(0) = I$, is the free product 
of certain groups ${\cal M}_z$ where $z$ is a $\nu$-vector
drawn from the set $\Xi
\stackrel{\rm def}{=} 
\{z = \mbox{column}\{\zeta_1,\zeta_2,...\zeta_{\nu} \} 
~|~ \zeta_1 = 1, z^*Dz = 0\}$.
Matrices in each ${\cal M}_z$ are explicitly
and uniquely parametrized in the paper.
Thus, every $\nu \times \nu$ matrix $U=U(\omega) \in {\cal M}_0$
can be 
represented in the form
$U = G_{z_1} \cdot G_{z_2} \cdot ... \cdot G_{z_\eta}$ 
with $\nu \times \nu$
polynomial matrix multiples 
$ G_{z_{\iota}} \in {\cal M}_{z_{\iota}}$, $z_{\iota} \in \Xi$,
for ${\iota}=1,...\eta$, so that
$z_{\iota} \ne z_{\iota+1}$ 
for ${\iota}=1,...\eta-1$
and this representation is unique.
The uniqueness 
includes the number of multiples $\eta$ where $\eta = 0, 1, 2, ... $,
their particular sequence 
$G_{z_1}, G_{z_2}, ... G_{z_\eta}$,
and the multiples themselves
with their respective parametrizations;
all these items can be defined in only one way
once the $U$ is given.

${ }$\\
${ }$\\
${ }$

{\bf AMS subject classification:} primary ??, secondary ??

\end{abstract}
% ----------------------------------------------------------------

\[
\begin{array}{c}
{ }
\end{array}
\]

{\bf 1.} Let 
\begin{equation}
\label{ueqsum}
U(\omega) = 
U_0
+ U_1\, \omega
+ U_2\, \omega^2
+ ...
+ U_{\kappa}\, \omega^{\kappa}
\end{equation}
be a polynomial matrix
(p.m.) of size $\nu_1 \times \nu_2$. 
Here $U_{\iota}$ for $\iota = 0,1,...\kappa$, 
are complex $\nu_1 \times \nu_2$ matrices,
$\omega$ is a scalar variable, and $\kappa = 0, 1, ...$.
Under assumption that variable $\omega$ 
takes on only real values,
$\omega^* = \omega$,
we extend the Hermitian conjugatation operation
\footnote{For a constant $\nu_1 \times \nu_2$ matrix $B$,
its Hermitean conjugate $B^*$ is defined
as the  $\nu_2 \times \nu_1$ transpose of the matrix
whose elements are complex conjugates of the elements
of $B$.
If $U_{\kappa} \neq 0$ in \eqref{ueqsum} then $\kappa$ is 
the degree of $U(\omega)$, denoted
$\kappa = \mbox{deg}\, U(\omega)$.
The following notations and conventions are assumed.
Greek lower case denotes scalars, Latin lower case denotes column vectors
of height $\nu$, e.g., 
$z = \mbox{column}\{\zeta_1,\zeta_2,...\zeta_{\nu} \}$.
Latin capitals denote matrices, $\mbox{i} = \sqrt{-1}$,
$I$ is the identity $\nu \times \nu$ matrix.
The symbol $\ni$ is read ``such that".
The notation $\{\Gamma_1~|~\Gamma_2 \}$ is used
for the set with elements $\Gamma_1$ satisfying the defining
property $\Gamma_2$.
}
to a p.m. $U(\omega)$
as follows:
\[
U(\omega)^* 
\stackrel{\rm def}{=}
U_0^*
+ U_1^*\, \omega
+ U_2^*\, \omega^2
+ ...
+ U_{\kappa}^*\, \omega^{\kappa}~.
\]
For $\nu \ge 2$, let $D$ be the diagonal 
$\nu \times \nu$ matrix of the form 
$D = \mbox{diag} \{-1,1,...1\}$.
In this paper we shall describe the group $\cal{M}$ of all p.m.'s
$U(\omega)$ of size $\nu \times \nu$ which satisfy the equation
\begin{equation}
\label{udu*=d}
  U(\omega) \cdot D \cdot U(\omega)^* = D,
\end{equation}
i.e., the group of p.m.'s unitary in the indefinite metric of index 1.
The description is effective in that every matrix
in $M$ is parametrized in a unique canonical form.
The main results are stated in \S \S 2, 3.
This problem relates to the problem \cite{Yakubovitch1}
of finding a factorization of a 
%non-positive definite 
p.m. $A(\omega) = A(\omega)^*$ 
in the form 
\begin{equation}
\label{a=ucu*}
A(\omega) = U(\omega) \cdot C \cdot U(\omega)^*
\end{equation}
with a p.m. $U(\omega)$ and a constant matrix $C = C^*$,
both of size $\nu \times \nu\,$.
Note that matrix $C$  in \eqref{a=ucu*} is not necessarily positive 
definite. 
As a special case, one might seek to factorize a constant matrix 
$A(\omega) = C = C^*\,,\mbox{det}\,C \neq 0$.
In this case, without loss of generality, 
we can assume $C$ 
to be a diagonal matrix 
\begin{equation}
\label{cpm1}
C = \mbox{diag} \{\pm 1, \pm 1,...,\pm 1 \}\,.
\end{equation}
The structure of the set of the  solutions $U(\omega)$ of the equation
\eqref{a=ucu*}, where $A(\omega) \equiv C$,
depends on how many $+1$'s and $-1$'s are 
on the diagonal of $C$.
Indeed, if the elements of the diagonal of $C$ are all of the same sign,
then equation \eqref{a=ucu*}, where $A(\omega) \equiv C$ ,
coupled with the normalizing condition $U(0) = I$,
has only the trivial solution $U(\omega) \equiv I$.
The set of the solutions becomes non-trivial when $C$ has both $+1$'s
and $-1$'s on the diagonal.
In this paper we consider the case of exactly one $-1$ with the rest
of diagonal elements being $+1$.
The same set of the solutions $U$ will result if $D$ is
substituted with $-D$, a diagonal matrix that has
exactly one $+1$ with the rest of diagonal elements being $-1$.
These two cases
exhaust the set of all indefinite matrices $C$
of the form \eqref{cpm1} in dimensions $\nu \le 3$.
Without loss of generality, we assume that the diagonal matrix $C$ has
the single $-1$ being at the top of the diagonal
and the rest of the diagonal consists of $+1$'s,
i.e., $C = D$ for the matrix $D$ defined above.
The task of finding a symmetric factorization \eqref{a=ucu*}
recurs in many contexts, for example in
the synthesis of linear optimal control
in differential games \cite{Yakubovitch2}.
\\

{\bf 2}. Introduce the sets:
\\
of vectors
\[
\Xi 
\stackrel{\rm def}{=}
\{z =
\mbox{column} \{\zeta_1 ,...,\zeta_{\nu}\}~|~ \zeta_1 = 1, z^*Dz = 0 \}~,
\]
\[
\Delta_z^0 
\stackrel{\rm def}{=}
\{
d~|~d^*z = d^* Dz = 0
\},~\mbox{for\ any\ } z \in \Xi\,,
\]
of vector polynomials
\[
\Delta_z 
\stackrel{\rm def}{=}
\{g = g_1 \omega + g_2 \omega^2 + ...~|~g_{\iota} \in \Delta_z^0 \},~
\mbox{for\ any\ } z \in \Xi\,,
\]
and of scalar polynomials
\[
\Phi
\stackrel{\rm def}{=}
\{
\phi = \phi(\omega)~|~\phi(0) = 0, \phi(\omega)^* = - \phi(\omega)
\}~.
\]

Then:
\\

1) the set of generatrices of the cone
\[
\Lambda \stackrel{\rm def}{=} \{x~|~x^*Dx = 0 \}
\]
is parametrized in a one-to-one way by $\Xi$ ;
\\

2) $\Delta_z \neq 0$ only if $\nu > 2$  and if 
$g \in \Delta_z$ then $g = Dg$;
\\

3) an element $\phi \in \Phi$ is of the form
$\phi(\omega) =
\displaystyle{
\mbox{i}\, \sum_{1\le \iota \le \kappa} \rho_{\iota} \omega^{\iota}
}
$,
where $\mbox{i} = \sqrt{-1}$,
the $\rho_{\iota}$ are real numbers, and $\kappa = 1, 2,...$.
\\

For $z \in \Xi$, $\phi \in \Phi$, $g \in \Delta_z$ define
the following  $\nu \times \nu$ polynomial matrix
\begin{equation}
\label{Gz}
G_z (\phi,g) 
\stackrel{\rm def}{=}
D [z(\phi - (1/2)\cdot g^*g)z^* + (zg^* - gz^*)] + I ~.
\end{equation}

The significance of the matrix in \eqref{Gz} is that, as will
be shown, any matrix $U \in \cal M$ can be decomposed
into p.m.'s of this form and a constant matrix.
\\

Define the sets of matrices:
\\

${\cal{M}}_z 
\stackrel{\rm def}{=} 
\{
G_z(\phi,g)~|~\phi \in \Phi, g \in \Delta_z
\}~ \mbox{for\ any\ } z \in \Xi$~,
\\

${\cal{M}}_0 
\stackrel{\rm def}{=} 
\{ 
U(\omega) \in {\cal{M}} ~ | ~ U(0) = I
\} 
$~,
\\

$
{\cal{N}}
\stackrel{\rm def}{=}
\{
V \in {\cal{M}} ~|~ \mbox{deg}\, V = 0
\}
$~,
\\

$
{\cal{\Upsilon}}
\stackrel{\rm def}{=}
\{
W \in {\cal{N}} ~ | ~ W = \mbox{diag}\{1,L\},
L \in \{L ~ | ~ L^* L = I_{(\nu-1) \times (\nu-1)}\} 
\}\,.
$
\\

Clearly, every p.m. $U(\omega) \in \cal{M}$ 
may be represented uniquely in the form
$U(\omega) = U_0 (\omega) \cdot V$ where 
$U_0 (\omega) \in {\cal{M}}_0 ,~ V\in \cal{N}$.
\\

The main result on the structure of $\cal{M}$ is the following
\\

{\bf Theorem 1.} 
{\em 1). For every}
$z \in \Xi$\,,
{\em the set}
${\cal {M}}_z$
{\em is a group under matrix multiplication.
The group}
${\cal M}_0$
{\em is the free product} 
\cite{Kuros}
{\em of groups }
${\cal {M}}_z$
{\em for all}
$z \in \Xi$.
{\em This means that any}
$U=U(\omega)\in{\cal M}_0$ 
{\em can be decomposed
in the form}
\begin{equation}
\label{udecomp}
U = G_{z_1} (\phi_1,g_1)\cdot ...
\cdot G_{z_{\eta}} (\phi_{\eta}\,,g_{\eta})~,
\end{equation}
{\em where} 
$G_{z_{\iota}} \in {\cal {M}}_{z_{\iota}}$ 
{\em for} 
$\iota = 1,...\eta$.
{\em Moreover,
if by aggregating 
the consecutive multiples
in the sequence in}
\eqref{udecomp}
{\em that belong to the same group
we make sure that}
$z_{\iota} \neq z_{\iota +1},~\iota = 1,...\eta-1$,
{\em then the obtained decomposition}
\eqref{udecomp}
{\em becomes unique for a given}
$U$.
{\em The uniqueness includes the number of multiples} 
$\eta$,
{\em the particular sequence}
$z_1\,,z_2\,, ... z_{\eta}$
{\em in} \eqref{udecomp}
{\em with the sequence}
${\cal {M}}_{z_1}\,,{\cal {M}}_{z_2}\,, ... {\cal {M}}_{z_{\eta}}$
{\em of the corresponding groups},
{\em and the multiples}
$G_{z_{\iota}} \in {\cal {M}}_{z_{\iota}} $ 
{\em themselves,
each of which is a unique p.m. for a given}
$U(\omega) \in {\cal M}_0$.
\\

{\em 2). Multiplication in each}
${\cal {M}}_z$, ~$z \in \Xi$,
{\em satisfies the condition}
\begin{equation}
\label{multinGz}
G_z (\phi,g) \cdot G_z (\psi, h) = 
G_z (\phi + \psi + (1/2)\cdot(h^*g - g^* h), g + h)\,,
\end{equation}
{\em where} $\phi,\psi \in \Phi$, ~ $g,h \in \Delta_z$.
\\

{\em 3). Different groups }
${\cal M}_z$ 
{\em are isomorphic; the connecting isomorphism is as follows}
\begin{equation}
\label{WGW}
W \cdot G_z (\phi,g) \cdot W^{-1} = G_{Wz} (\phi, Wg) \,,
\end{equation}
{\em where} 
$W \in {\cal \Upsilon},~ z \in \Xi,~\phi\in \Phi,~g\in \Delta_z$,
{\em and necessarily}
$W^{-1} \in {\cal \Upsilon},~ Wz \in \Xi$
{\em and}
$Wg\in \Delta_{Wz}$.
\\

{\em 4). The mapping}
$(z,\phi,g) \mapsto G_z (\phi,g)$
{\em is one-to-one on the set}
\[
\{
(z,\phi,g)~|~z \in \Xi,~\phi\in \Phi,~g\in \Delta_z,~\phi \neq 0~ \mbox{or} ~ g \neq 0
\}
\]
{\em If}
$U(\omega)$
{\em is decomposed in the form}
\eqref{udecomp}
{\em with}
$z_{\iota} \neq z_{\iota +1},~\iota = 1,...\eta-1$,
{\em then}
\[
\mbox{deg}\, U(\omega) = 
\mbox{deg}\,G_{z_1} (\phi_1,g_1) + ... + 
\mbox{deg}\, G_{z_{\eta}} (\phi_{\eta},g_{\eta})
= \mbox{deg} (\phi_1 - g_1^*g_1) + ... 
+ \mbox{deg} (\phi_{\eta} - g_{\eta}^*g_{\eta})\,.
\]
{\em Here} 
$z_{\iota} \in \Xi,~\phi_{\iota} \in \Phi, 
~ g_{\iota} \in \Delta_{z_{\iota}} ,~\iota = 1,...\eta$.
\\

{\em 5). The center of}
${\cal M}_z , ~z \in \Xi$,
{\em is the set} 
$\Psi_z 
\stackrel{\rm def}{=} 
\{G_z (\phi,0)~|~\phi \in \Phi
\}$
{\em If}
$\nu > 2$,
{\em then the commutant of}
${\cal M}_z$
{\em coincides with the center of}
${\cal M}_z$.
{\em If}
$\nu = 2$,
{\em then the group }
${\cal M}_z$
{\em is commutative and coincides with}
$\Psi_z$\,.
\\

{\bf 3. Real cases.}
Let $\cal M'$ be the subgroup of $\cal M$ which consists of p.m.'s with
real coefficients.
Let $\Xi' ,~\Phi' ,~\Delta'_z ,~z\in \Xi'$ be the real analogs of the
sets $\Xi ,~\Phi ,~{\Delta}_z$ introduced in \S 2.
Clearly, $\Phi' = \{0\}$.
As in \S 2 , let us define ${\cal M'}_z$ and ${\cal M}_0$.
The structure of the group ${\cal M}_0$ is described in the following
\\

{\bf Theorem 2.} 
{\em If} $\nu = 2$, 
{\em then the group} 
${\cal M}_0$ 
{\em is trivial. If}
$\nu > 2$,
{\em then the group}
${\cal M}_0$
{\em is not trivial and is the free product of groups}
${\cal M'}_z$
{\em for all}
$z\in \Xi'$.
{\em Groups}
${\cal M'}_z\,,~z\in \Xi'$,
{\em are commutative}.
\\

It is also interesting to describe the subgroup of p.m.'s in $\cal M$
with real coefficients not with respect to variable $\omega$,
but with respect to variable $\lambda = \mbox{i} \omega$.
Note, that whereas $\omega^* = \omega$,
we have $\lambda^* = -\lambda$.
Thus, a p.m. \eqref{ueqsum} has real coefficients
with respect to variable $\omega$ if all $U_{\iota}$ are real
matrices; a p.m. \eqref{ueqsum} 
has real coefficients with respect to variable $\lambda$
if  $U_{\iota}$ are real for even ${\iota}$ 
and $\mbox{i} U_{\iota}$ are real for odd ${\iota}$.
As above, we introduce 
$\Phi"=\{\phi = \phi(\lambda)~|~\phi(0) = 0, \phi^* = -\phi\}$.
The general form of such a $\phi$
is 
$\phi(\lambda) = \rho_1 \lambda + \rho_2 \lambda^3 + \rho_3 \lambda^5 + ...$,
where $\rho_{\iota}$ are real numbers, and the sum is finite.
In the $\omega$ representation the general form
of a $\phi \in \Phi"$ is ~
i$(\rho_1 \omega + \rho_2 \omega^3 + \rho_3 \omega^5 + ...)$,
with real $\rho_{\iota}$.
We define $\Delta_z^" ,~{\cal M}_z^"$ for $z \in \Xi"$ and also ${\cal M}_0^"$.
\\

{\bf Theorem 3.}
{\em 1). The group}
${\cal M}_0^"$
{\em is a free product of groups}
${\cal M}_z^"$
{\em for all}
$z \in \Xi"$.
\\
{\em 2). If}
$\nu = 2$,
{\em then}
$\Xi"$ 
{\em consists of two elements}
$\Xi" = \{z_1 ,~z_2\}$, 
$z_1 = \mbox{column}\,\{1,1\}$
{\em and}
$z_2 = \mbox{column}\,\{1,-1\}$,
{\em so that the free product just mentioned contains two groups}
${\cal M}_z^"$
{\em only, and each of these two groups is commutative};
\\
{\em 3). If}
$\nu > 2$,
{\em then each group}
${\cal M}_z^"$
{\em is not commutative. The commutant of}
${\cal M}_z^"$
{\em coincides
with the center and is equal to}
$\{G_z (\phi,\,0)~|~\phi \in \Phi"\}$.
\\
  
As an example of application of Theorem\ 3, let us
present a general form of a $2 \times 2$ p.m. $U(\lambda) \in {\cal M}_0^"$
of degree at most 2.
We have: either 
$U(\lambda) = G_{z_1} (\alpha \lambda) \cdot G_{z_2} (\beta \lambda)$
or
$U(\lambda) = G_{z_2} (\alpha \lambda) \cdot G_{z_1} (\beta \lambda)$.
Here $z_1$ and $z_2$
are the two elements in $\Xi'$ when $\nu = 2$;
$\alpha \lambda$ and $\beta \lambda$ are $\lambda$-monomials
in $\Phi^"$ with real coefficients $\alpha$ and $\beta$. 
Each monomial
is of degree 1, if the coefficient is non-zero. 
Otherwise
the monomial is zero.
We can go further and represent these two types of $U(\lambda)$
in the per-component form.
In the first case, we calculate that
\[
U(\lambda) = \left [
\begin{array}{l l}
-2\alpha\beta\lambda^2  - (\alpha + \beta)\lambda + 1~~~ &
~~2\alpha\beta\lambda^2  - (\alpha - \beta)\lambda \\
~~2\alpha\beta\lambda^2 + (\alpha - \beta)\lambda~~~ &
-2\alpha\beta\lambda^2  + (\alpha + \beta)\lambda + 1 
\end{array}
\right ]\,,
\]
and in the second case,
\[
U(\lambda) = \left [
\begin{array}{l l}
-2\alpha\beta\lambda^2  - (\alpha + \beta)\lambda + 1~~~ &
-2\alpha\beta\lambda^2  + (\alpha - \beta)\lambda \\
-2\alpha\beta\lambda^2 - (\alpha - \beta)\lambda~~~ &
-2\alpha\beta\lambda^2  + (\alpha + \beta)\lambda + 1 
\end{array}
\right ]\, .
\]
\\

{\bf Proofs.} Let $a\|b$ denote the existence of a linear dependence 
between vectors $a$ and $b$;
we shall write $a\|\!|b$ to mean that $a\|b$ and,
moreover, that the coefficients of the nontrivial linear combination
of $a$ and $b$ can be chosen real.
The proofs of Lemmas 1 to 2 below present no difficulties.
\\

{\bf Lemma 1.} $(a,b \in \Lambda)
~\Longrightarrow~
((a^*Db = 0) \Longleftrightarrow (a\|b))\,.$
\\

{\bf Lemma 2.} $(ab^* = ba^*)
\Longleftrightarrow
(a\|\!|b)\,.$
\\

{\bf Lemma 3.} $(XDX^* = X^*DX = 0)
~\Longrightarrow~
(\exists y, z \in \Xi, \alpha \ni X = \alpha D yx^*)\,.$
\\

The lemma is proved by applying Lemma 1 
to the rows and, separately, columns of $X$.
\\

{\bf Lemma 4.} 
{\em Let}
$y,z \in \Xi$
{\em and a)}
$XDz = 0$;
{\em b)}
$\forall k \in \Delta_z^0 ~
(XDk \| Dy)$.
{\em Then}
$\exists s \in \Delta_z^0 ,~ r \in \Delta_y^0, ~ \alpha \ni
X = D(rz^* - ys^* + \alpha yz^*)\,.$
\\

{\bf Lemma 5.}
{\em Let}
$z_{\iota} \in \Xi\,,~\iota = 1,...,\eta$.
{\em Then}
\[
(Dz_1z_1^* \cdot Dz_2z_2^* \cdot ... \cdot Dz_{\eta}z_{\eta}^* = 0)
\Longleftrightarrow
(\exists 
\iota_0 
(1\le \iota_0 < \eta ) 
\ni 
z_{\iota_0} 
= z_{\iota_0 + 1} 
)\,.
\]

The product of diadic matrices $Dz_{\iota}z_{\iota}^*$, $\iota = 1,...\eta$,
in the lemma is equal to $(Dz_1 z_{\eta}^*)\xi$,
where 
$\xi = (z_1^*Dz_2) \cdot ... \cdot (z_{\eta -1}^* Dz_{\eta} )$.
But $Dz_1 z_{\eta}^* \neq 0$ and the condition $\xi = 0$
is equivalent to the existence of an $\iota$,
$1 \le \iota \le \eta -1$,
such that $z_{\iota}^* D z_{\iota + 1} = 0$.
Since $z_{\iota} , z_{\iota + 1} \in \Lambda$,
it follows from Lemma 1 that $z_{\iota} \| z_{\iota + 1}$;
hence, in view of the normalization, $z_{\iota} = z_{\iota + 1}$,
and $\Longrightarrow$ implication is true.
Obviously the 
$\Longleftarrow$
implication is also true.

{\bf Lemma 6.}
{\em Let}
$\phi \in \Phi,~g \in \Delta_z, z \in \Xi$
{\em and suppose}
$\phi$
{\em and} $g$ 
{\em are not both zeros.}
\\
{\em Then 1) the leading coefficient of a p.m.}
$G_z (\phi ,g)$
{\em is proportional to}
$Dzz^*$;
{\em 2)}
$\mbox{deg}\,G_z(\phi ,g) = \mbox{deg}\,(\phi - g^*g) > 0$.
\\

The most laborous is the proof of
\\

{\bf Lemma 7.}
{\em If}
$\mbox{deg}\,U(\omega) > 0$,
{\em then the degree of the p.m.}
$U(\omega) \in \cal M$
{\em may be decreased by a right or left multiplication
by a p.m. of the form}
\eqref{Gz}
{\em (and real if}
$U(\omega)$ 
{\em is real).}
\\

{\bf Proof of Lemma 7.}
Let
$U(\omega) = \sum_{\iota=0}^{\kappa} X_{\kappa - \iota}\, \omega^i$\,,
where $X_{\iota}\,,\iota = 0,1,...\kappa$, are complex
(or real, or $\lambda$-real, depending on the case) constant $\nu \times \nu$
matrices, and $X_0 \ne 0$, so that
$\mbox{deg}\, U(\omega) = \kappa$.
It follows from
\eqref{udu*=d}
that $U(\omega)^* \cdot D \cdot U(\omega) = D$.
Extending definition of $X_{\iota}$ for
$\iota > \kappa$ to be null matrices,
we obtain a family of equalities
\[
~~~~~~~~~~~~~~~~~~~~~~~~~~~~~
X_0\,D\,X_{\gamma}^* + X_1\,D\,X_{\gamma-1}^* + ...
+ X_{\gamma}\,D\,X_0^*  = 0\,
~~~~~~~~~~~~~~~~~~~~~~~~~~~~~(9_{\gamma})
\]
\[
~~~~~~~~~~~~~~~~~~~~~~~~~~~~~
X_0^*\,D\,X_{\gamma} + X_1^*\,D\,X_{\gamma-1} + ...
+ X_{\gamma}^*\,D\,X_0  = 0\,
~~~~~~~~~~~~~~~~~~~~~~~~~~~~~(10_{\gamma})
\]
Here $\gamma = 0, 1, ..., 2\kappa -1$.
\addtocounter{equation}{+2}
Having ($9_0$) and ($10_0$) we can apply Lemma 3 to $X = X_0$.
For the leading coefficient $X_0$,
this yields its diadic representation 
$X_0 = \alpha_0 Dyz^*$ for some vectors
$y ,z \in \Xi\,,$ and a number $\alpha_0$\,.
Note that $\alpha_0 \neq 0$\, 
because $\mbox{deg}\, U(\omega) = \kappa$ .
The $y$ and $z$ will be our candidates for the subscript 
of the p.m. of the form \eqref{Gz} which should
decrease the degree of $U(\omega)$ after multiplying
$U(\omega)$ on the left or on the right, respectively.
Specifically, we will prove the lemma
if we show that at least one of the two possibilities holds:
\\
I) $\exists \phi \in \Phi, \exists g \in \Delta_z \ni \mbox{deg}\, 
[U(\omega)\cdot G_z(\phi,g)] < \kappa$, 
\\
II) $\exists \psi \in \Phi, \exists h \in \Delta_y \ni \mbox{deg}\, 
[G_y(\psi,h)\cdot U(\omega)] < \kappa$\,.
\\
(In these statements, sets $\Xi, \Phi, \Delta_x, \Delta_y$
should be appropriately  modified in the cases of reals.)
\\

For integer positive numbers $\tau,\mu$, and $\xi$,
such that $\tau \le \mu$ and $\tau \le  \xi$, 
consider the following conditions:
\\
A)$~\exists \alpha_{\iota}\,, 0 \le \iota \le \tau -1\,, 
\ni X_{\iota} = \alpha_{\iota} D y z^* \,$,
i.e., the first $\tau$ leading coefficients of p.m. $U(\omega)$
are diadic matrices proportional to $X_0$;
\\
B)$~X_{\iota}Dz = 0 ~\mbox{for\ all}~ \iota, \tau \le \iota \le \mu-1,~ 
\mbox{and}~X_{\mu} Dz \neq 0\,$,
i.e., 
the 
$\mu -\tau$ coefficients of p.m. $U(\omega)$,
that follow the last coefficient $X_{\tau -1}$ mentioned in A),
turn into 0, 
when multiplied on the right by the vector-column $Dz$, 
but this does not haappen for the $\mu -\tau +1$st coefficient;
\\
C)$~y^*X_{\iota} = 0, ~\mbox{for\ all}~ \iota, \tau \le \iota \le \xi-1,~
\mbox{and}~y^* X_{\xi} \neq 0\,$,
i.e., 
the 
$\xi -\tau$ coefficients of p.m. $U(\omega)$,
that follow the last coefficient $X_{\tau -1}$ mentioned in A),
turn into 0,
when multiplied on the left by the vector-row $y^*$,
but this does not happen for the $\xi -\tau +1$st coefficient;
\\
D)$~\exists\, s \in \Delta_z^0\,,r \in \Delta_y^0\,,\alpha_{\tau} \ni
X_{\tau} = D(rz^* - ys^* + \alpha_{\tau} yz^* )\,$,
i.e., the coefficient $X_{\tau}$ is the sum of three diadic matrices
as stated.
\\

Consider the largest $\tau \ni$ A). 
Since $X_0 = \alpha_0 D y z^*$, the $\tau$ is at least 1.
It can not be larger than $\kappa$, though.
Hence, $1 \le \tau \le \kappa$.
Let $\mu \ni$ B), $\xi \ni$ C).
Clearly, $\tau \le \mu,\xi \le \kappa$.
The way the proof proceeds further depends of whether the two inequalities
\begin{equation}
\label{mugt2tau}
\mu \ge 2 \tau \mbox{ and }  \xi \ge 2 \tau 
\end{equation}
both hold or not.
First we consider 
the easier 
\\

{\em Case 1: at least one 
inequality in } \eqref{mugt2tau} {\em fails}.
Suppose, for example, that 
%the first inequality in
%\eqref{mugt2tau} fails, that is, we have 
$\mu < 2\tau$.
We will then verify I).
Equality ($9_{\mu}$) and Lemma 2 imply
that there exists a real number
$\rho$ such that
\begin{equation}
\label{cancel}
\alpha_0 D y + \mbox{i} \rho X_{\mu} Dz = 0\,.
\end{equation}
Using \eqref{cancel} it is easy to verify
that 
$\mbox{deg}\,[U(\omega) \cdot (Dz(\mbox{i}\rho \omega^{\mu} ) z^* + I) ] < \kappa$,
i.e., as stipulated in I), the degree decreases when
$U(\omega)$ is multiplied on the right by 
$G_z(\phi,g)$ with $\phi(\omega)$ and $g(\omega)$ taken here
as
 $\phi(\omega) = i\rho\omega^{\mu}\,, g(\omega)=0$.
(When $z, y$ and $X_{\mu}$ are real, 
since $\rho$ is real, \eqref{cancel} implies 
that i$\alpha_0$ is real and hence
all elements of i$X_0$ are real. It follows, 
that if $U(\omega)$ has only real coefficients,
inequality $\mu < 2\tau$ can not occur.
In case of reals with variable $\lambda = i \omega$,
inequality $\mu < 2\tau$ can only occur
for odd $\mu$.)
Analogously
we verify II) if the second inequality in \eqref{mugt2tau}
fails.
\\

{\em Case 2: both
inequalities in}
\eqref{mugt2tau} 
{\em hold}.
This case will be more laborous.
First,
we will verify D) by using Lemma 4 where $X = X_{\tau}$\,.
The conditions a) and b) in Lemma 4 are obviously satisfied.
Let us verify the condition c).
Pick a vector $k \in \Delta_z^0$ and denote 
$v \stackrel{\rm def}{=} X_{\tau} Dk$. 
Condition c) is then rewritten as
$v \| Dy$.
To verify the latter,
we are going to use Lemma 1 with $a = Dy$, $b = v$.
Obviously $a \in \Lambda$, so we have to check that also
I') $b \in \Lambda$, i.e., $v^* Dv = 0$, and that
II') $a^*Db = 0$, i.e., $y^*v = 0$.
\\

Multiply the equality ($9_{2\tau}$) by $Dk$ on the right
and by $k^*D$ on the left.
It follows from \eqref{mugt2tau} that
on the left-hand side all the summands but one will turn to zero.
In other words, I') holds.
We also have
$y^*v = (y^* X_{\tau} ) D k = 0$,
that is, II') holds.
\\

Now, as the pre-conditions of Lemma 4 are satisfied,
we obtain
representation 
\begin{equation}
\label{xtau}
X_{\tau} = D(rz^* - ys^* + \alpha yz^*)\,.
\end{equation}
Observe, that vectors $s$ and $r$ can not
be both zero in \eqref{xtau},
because it would have contradicted
to the definition of $\tau$
which could have been possible to increase in such a case.
Let us assume that $s \neq 0$ in \eqref{xtau}
and show I).
If $r \neq 0$ in \eqref{xtau} we can similarly
show II).
\\

Define
$w
\stackrel{\rm def}{=}
X_{2\tau} D z\,,
p
\stackrel{\rm def}{=}
\alpha_0 D y$.
Using Lemma 1, we will show that
$w\|p$.
Obviously $p^*Dp = 0$ and
we will also show that
I") $w^* Dw = 0$
II") $p^*Dw = 0$
which are the preconditions in Lemma 1
applied to $a = w$, $b = p$.
\\

Using ($10_{4\tau}$) immediately yields I").
To obtain II"), we begin with ($10_{2\tau}$)
which we multiply 
on the right
at $Dz$.
This yields $\alpha_0^* z y^* w = 0$
which implies II").
Now we can apply Lemma 1 and obtain
\begin{equation}
\label{wp}
w = (\sigma_0 + \mbox{i} \rho_0 ) \cdot p \,,
\end{equation}
where $\sigma_0$ and $\rho_0$ some real numbers.
(In the case of reals with respect to variable $\omega$
the resulting $\rho_0$ will be zero,
in the case of reals with respect to variable $\lambda$
the resulting $\rho_0$ will be zero for ?????).
\\

We can now use the obtained expressions for
substituting in ($9_{2\tau}$) 
\[
X_0 = pz^*,~ 
X_{2\tau} = (\sigma_0 + \mbox{i} \rho_0 ) \cdot p\,~
X_{\tau} = -ps_0 + nz^*
\] 
where we denoted 
$n
\stackrel{\rm def}{=}
D(r + \alpha_{\tau} y)$
and
$s_0 
\stackrel{\rm def}{=}
s/\alpha_0^*$.
We then obtain an equation
$pp^* (2 \sigma_0 + s_0^* D s_0 ) = 0$
from which,
taking into account equality
$s_0^* D s_0 = | s_0 |^2$,
we obtain $\sigma_0 = - (1/2) | s_0 |^2$.
\\

We are now ready to verify I).
Set 
$\phi (\omega) 
\stackrel{\rm def}{=}
\mbox{i} \rho \omega^{2\tau}\,$
and
$g(\omega) 
\stackrel{\rm def}{=}
d \omega^{\tau}$,
where $d \in \Delta_z^0$ and real $\rho$ are
constants to be defined.
The definition should satisfy
condition
\begin{equation}
\label{toeps}
X_0 +
X_{2\tau}Dzz^* (\mbox{i} \rho - (1/2) |d|^2)
+X_{\tau} D (zd^* -dz^* ) 
\stackrel{\rm def}{=}
pz^* \epsilon = 0  .
\end{equation}
Parameter $\epsilon$ which is defined in
\eqref{toeps} 
is equal to
\[
\epsilon 
=
1 + s_0^* d + 
(\mbox{i} \rho_0 - (1/2) |s_0|^2 )
\cdot
(\mbox{i} \rho - (1/2) |d|^2 ).
\]
Our goal is to set parameters $\rho$ and $d$ so that
$\epsilon$ will become zero.
We are looking for $d$ of the form
$d = \theta s_0$,
where $\theta_0$ is a constant to be defined.
It follows, that
$\epsilon = \epsilon_1 - \mbox{i} \cdot (1/2) |s_0|^2 \cdot \epsilon_2$.
where
$\epsilon_1 = |1 + (1/2) \theta \cdot |s_0|^2 |^2 - \rho \cdot \rho_0$
and
$\epsilon_2 = \rho_0 \cdot |\theta|^2  - 2 \cdot \mbox{Im} \theta + \rho$.
If $\rho_0 = 0$, we set $\rho = 0$ and $\theta = -2|s_0|^{-2}$
to obtain $\epsilon = 0$.
If $\rho_0 \neq 0$, then by setting
$\rho = |1 + (1/2) \theta \cdot |s_0|^2|^2 /\rho_0$
we obtain $\epsilon_1 = 0$.
Such setting of $\rho$ results in
\begin{equation}
\label{eps2}
\epsilon_2 = \rho_0 ^{-1} \{
1 +
|\theta |^2 
\cdot 
[ \rho_0^2 + ((1/2) |s_0|^2)^2] 
-2 [ \rho_0 \cdot \mbox{Im} \theta - 
((1/2)|s|^2) \cdot \mbox{Re} \theta] 
\}
\end{equation}
Now setting
\[
\theta = 
(\mbox{i} \rho_0 - (1/2) |s_0|^2 )
\cdot 
(\rho_0^2 + ((1/2) |a_0 |^2)^2)^{-1}
\]
results in $\epsilon_2 = 0$ as required.
Lemma 5 is proved.
\\

{\bf Proof of Theorem 1.} 
Let us first establish \eqref{multinGz}.
We have
\[
G_z (\phi,g) \cdot G_z (\psi, h) = (D A + I) (D B + I) =
DADB + D(A+B) + I\,,
\]
where $A = z(\phi - (1/2)\cdot g^*g)z^* + (zg^* - gz^*)$
and
$B = z(\psi - (1/2)\cdot h^*h)z^* + (zh^* - hz^*)$
according to \eqref{Gz}.
Using identities $z^*Dz = 0$, $z^*Dh = 0$, and $g^*Dz = 0$,
that follow from the definitions, and multiplying
each summand in $A$ by each summand in $B$ on the right,
we establish that the only 
non-zero term in the resulting sum for $DADB$ is 
$-Dzg^*Dhz^* = -Dz(g^*h)z^*$ where $Dh$ is substituted with $h$
since $h \in \Delta_z$.
Using the identity,
\[(g+h)^* (g+h) - (h^*g - g^*h) = g^*g + h^*h + 2g^*h\,,\] 
the sought product 
then becomes
\[
G_z (\phi,g) \cdot G_z (\psi, h) = -Dzg^*hz^*
+ D(A+B) + I = ~~~~~~~~~~~~~~~~~~~~~~~~~~~~~~~~~
\]
\[
= D [ z(\phi+\psi - (1/2)\cdot (g^*g + h^*h - 2g^*h))z^* + (z(g+h)^* - (g+h)z^*)
\]
\[
= G_z (\phi+\psi + (1/2)\cdot (h^*g - g^*h),g+h)
\]
as required in \eqref{multinGz}.
Similar and even simpler arguments show that for $U = G_z (\phi,g)$
\eqref{udu*=d} holds.
Hence sets ${\cal{M}}_z$ for $z \in \Xi$
are indeed groups under the matrix multiplication
and ${\cal{M}}_z \subset {\cal{M}}_0$ for any $z \in \Xi$.
\\

%Why each Gz is uniquely parametrized in that form?
It is easy to see that 
$\mbox{deg}\,G_z (\phi,g) = \mbox{deg} (\phi  - g^*g)$
unless both $\phi$ and $g$ are identical zeros in which case
$G_z (\phi,g)$ is identical to the matrix $I$.
Hence if $G_z (\phi,g)$ is not identical to $I$, then
the leading coefficient of $G_z (\phi,g)$ is proportional
to the diadic matrix $Dzz^*$.
If a product of matrices $G_z (\phi,g)$ is formed,
like that in the right-hand side of \eqref{udecomp}
with adjacent multiples belonging to different
groups ${\cal{M}}_z$, then the leading term of this product
will be proportional to the product of the diadic
$Dzz_{\iota}^*$ as in the left-hand side of the
statement in Lemma 5.
According to that lemma, the leading term can not degenerate
to zero. Because if it did, then a pair of adjacent multiples 
would have belonged to the same group ${\cal{M}}_z$.
Therefore, the products like those in the right-hand side of \eqref{udecomp}
can never degenerate to the identity matrix $I$,
given that the number of multiples is not zero and
that adjacent multiples belong to different groups ${\cal{M}}_z$.
Hence the product of group ${\cal{M}}_z$ where
$z$ runs over set $\Xi$, we temporarily
denote this product as $\tilde{\cal{M}}$, is {\em free}
(i.e., has no relations, see, e.g., \cite{Kuros} for the definition).
Lemma 7 tells that ${\cal{M}}_0 \subset \tilde{\cal{M}}$ and
obviously we have $\tilde{\cal{M}} \subset {\cal{M}}_0$,
thus $\tilde{\cal{M}} = {\cal{M}}_0$ and statement 1) 
of theorem 1 is proved.

{\bf  Acknowledgment}.
The author is indebted to V. A. Yakubovitch and D.K. Faddeev
for their interest in this paper.
% ----------------------------------------------------------------
%\newpage

% ----------------------------------------------------------------
\end{document}